\documentclass[12pt]{amsart}
\usepackage{amsmath}
\usepackage{amsfonts, amssymb, amscd}
\usepackage{xypic}

\def\ZZ{{\mathbb Z}}

\def\PP{{\mathbb P}}
\def\CC{{\mathbb C}}

\def\cF{\mathcal{F}}
\def\cO{\mathcal{O}}
\addtocounter{footnote}{1}

\newtheorem{lemma}{Lemma}[section]
\newtheorem{theorem}[lemma]{Theorem}
\newtheorem{corollary}[lemma]{Corollary}

\newtheorem{prop}[lemma]{Proposition}
\theoremstyle{definition}
\newtheorem{definition}[lemma]{Definition}

\newtheorem{remark}[lemma]{Remark}

\theoremstyle{remark}
\newtheorem*{proof*}{Proof}
\numberwithin{equation}{section}

\title[Free actions on c.i. of quadrics]{Classification of free actions on complete intersections of four quadrics}

\author{Zheng Hua}
\address{Department of Mathematics \\ University of Wisconsin \\
  Madison \\ WI \\ 53706 \\ USA\\
{\tt hua@math.wisc.edu}  }

\begin{document}

\begin{abstract}
In this paper we classify all free actions of finite groups on
Calabi-Yau complete intersection of 4 quadrics in $\PP^7$, up to
projective equivalence. We get some examples of smooth Calabi-Yau
threefolds with large nonabelian fundamental groups. We
also observe the relation between some of these examples and moduli of
polarized abelian surfaces.
\end{abstract}

\maketitle
\section{Introduction}
The original motivation of this paper is to generalize
Beauville's construction of Calabi-Yau manifolds with a non-abelian
fundamental group\cite{B}. As one result of this paper, we
construct many new examples of Calabi-Yau manifolds with
non-abelian fundamental groups. In particular we construct five families of Calabi-Yau threefolds
with fundamental groups of order 64. All these families are related to pencils of certain abelian surfaces.
Three of these families have been previously studied in \cite{GPa} and \cite{BH}.
The new examples are constructed as free quotients of small resolutions of singular complete intersections
of four quadrics in $\PP^7$ that contain a pencil of (2,4) polarized abelian surfaces (theorem \ref{24}).

We also classify all families of complete intersections
of four quadrics in $\PP^7$ with a free finite group action and at most ODP singularities.
The key idea is to use Holomorphic Lefschetz formula to obtain restriction on possible group actions.
This paper is quite elementary, the reasoning is sometimes very explicit and is
never very deep. Calculations of this paper can be generalized to other complete intersections in projective spaces
or in products of projective spaces.

The paper is organized as follows. In section \ref{Beauville}, we
review the construction of a smooth Calabi-Yau threefold with
quaternion group $H_8$ acting freely on it due to Beauville. We will
see how the character theory of $H_8$ and holomorphic Lefschetz
formula make this the only possible family of complete intersections
with $H_8$ action. We also see that no linear action of the dihedral
group $D_8$ could lead to any similar examples. In section
\ref{projrep}, we give a brief review about projective
representations of finite groups and define the terminology of
allowable actions, semi-allowable actions and Lefschetz condition.
Section \ref{alg} contains a scheme of the algorithm of classifying
(semi-)allowable actions on complete intersections of four quadrics
in $\PP^7$. As an application we make several tables in the next
section, listing all the (semi-)allowable actions with groups of
order from 2 to 64. In section \ref{ci} we compute the cut out
equations of families of Calabi-Yau threefolds with order 64
semi-allowable actions. There are two such families with five
different order 64 semi-allowable actions. In the last section we
prove the existence of equivariant small resolutions (\ref{p2} and
\ref{p11}). We also explain the relations between these Calabi-Yau
threefolds and moduli of polarized abelian surfaces.

All the group theoretic calculations are done in GAP\cite{GAP}. The software package Macaulay 2[M]
is also very useful to us in checking smoothness. I am grateful to my advisor Lev
Borisov, who gave many important ideas for this project.

\section{Beauville's example}\label{Beauville}
In this section we will first review Beauville's example of a free
action of quaternion group $H_8$ on a nine dimensional family of
smooth complete intersections of four quadrics in $\PP^7$(See
\cite{B}). Additionally we will explain why there is no such family
with free action of the dihedral group $D_8$. In the process we will
see how holomorphic Lefschetz formula leads to restriction on
possible free group actions.

The quaternion group $H_8$ is the group of order 8 with elements
${\pm1,\pm i,\pm j, \pm k }$ and $i^2=j^2=k^2=-1,ij=k,jk=i,ki=j$. By
a character calculation, $H_8$ has 4 one dimensional irreducible
representations and 1 two dimensional irreducible representation. We
denote them by $V_1,\ldots,V_4$ and $W$. The regular representation
V has decomposition $V=V_1\oplus V_2\oplus V_3\oplus V_4\oplus
W^{\oplus 2}$. The induced representation on the second symmetric
product of V has decomposition $Sym^2(V)=V_1^{\oplus 5}\oplus
V_2^{\oplus 5}\oplus V_3^{\oplus 5}\oplus V_4^{\oplus 5}\oplus
W^{\oplus 8}$. Pick 4 generic quadrics $q_1,\ldots,q_4$ such that
$q_i$ belongs to $V_i$. For generic choice of $q_i$, Beauville
showed that the complete intersection X in $\PP(V^*)$, given by
$q_1=\ldots=q_4=0$ is smooth and action of $H_8$ on X has no fixed
points. As a consequence the quotient variety $X/{H_8}$ is a smooth
Calabi-Yau manifold with fundamental group $H_8$.

The following theorem is a special case of the standard holomorphic
Lefschetz formula:
\begin{theorem}
Let $X$ be a smooth algebraic variety over $\CC$ and $f:X\rightarrow
X$ be a holomorphic automorphism of finite order with no fixed
points. For a linearized coherent sheaf $\cF$, the Lefschetz number
$$\Lambda(f,\cF)\colon=\sum^{m}_{q=0} (-1)^q Tr(f^{*};H^q(X,\cF))$$
is zero, where Tr stands for the trace.
\end{theorem}

Holomorphic Lefschetz formula explains why Beauville needed to pick
this particular representation $V$ and these particular choices of
quadrics $q_i$. We identify the vector space $V$ with
$H^0(X,\cO(1))$. By Kodaira's vanishing theorem and holomorphic
Lefschetz formula, $Tr(g,H^0(X,\cO(1)))$=0 for any $g$ non-identity.
The quaternion group $H_8$ has five conjugacy classes represented by
$\{(1),(i),(j),(-1),(k)\}$. By computing traces of each conjugacy
class, we get the trace vector $[8,0,0,0,0]$ for $V$, which means it
must be the regular representation. The induced representation
$Sym^2(V)$ has trace vector $[36,0,0,4,0]$. By Lefschetz formula
$H^0(X,\cO(2))$ has trace vector $[32,0,0,0,0]$. Their difference
$[4,0,0,4,0]$ is the trace vector for the space of four quadrics.
This is an actual group character for $H_8$. More precisely,
$[4,0,0,4,0]$ is the sum of characters of the $4$ one dimensional
irreducible representations $V_1,\ldots,V_4$. This is why Beauville
picked $q_i$ from the direct sum of copies of $V_i$ in $Sym^2(V)$.

The only other non abelian group of order 8 is the dihedral group
$D_8$. It's natural to ask that whether $D_8$ acts freely on any
smooth complete intersections of four quadrics in $\PP^7$. Dihedral
group $D_8$ is presented by $\{a,b|a^4=1,b^2=1;ab=ba^3\}$. It has
five conjugacy classes $\{(1),(b),(ab),(a^2),(a)\}$. Again, we
identify $V$ with $H^0(X,\cO(1))$, and we assume $\cO(1)$ can be
linearized so that $D_8$ acts on $V$. If $D_8$ acts freely on $X$,
the trace vector of $V$ should be $[8,0,0,0,0]$, i.e. $V$ must be
the regular representation. The trace vector for $Sym^2(V)$ is then
$[36,4,4,4,0]$. Subtracting $[32,0,0,0,0]$, we get $[4,4,4,4,0]$. It
is $\emph{not}$ a character of $D_8$. So $D_8$ can not act linearly
on any smooth complete intersection of four quadrics in $\PP^7$.

For any group $G$ of order bigger than eight, $\cO(1)$ can't be
$G$-linearized. Because otherwise the holomorphic Lefschetz formula
shows that the character of the action on $V=H^0(X,{\mathcal O}(1))$
is a fractional multiple of the character of the regular
representation, which leads to a contradiction. Hence instead of
linear representations we should look for projective
representations. In next section, we will give a brief review on
projective representations of finite groups. We will see how
holomorphic Lefschetz formula puts restriction on these projective
representations.

\section{Preliminaries of Projective Representations}\label{projrep}
In the first part of this section we recall some facts about
projective representations of finite groups. Our notations
follow\cite{Ber}. After that we define the notion of \emph{allowable
action} of a subgroup of $\mathbb{PGL}(8,\CC)$.

\begin{definition}
Let $G$ be a finite group. A triple $(\Gamma,f,A)$ is called a
central extension of $G$ if $\Gamma$ is a group, $A\subseteq
Z(\Gamma)$ and $f$ is a homomorphism of $\Gamma$ onto $G$ such that
$ker f = A$. A central extension $(\Gamma,f,A)$ is called Schur
Cover of $G$ if A equals the second group homology $H_2(G,\ZZ)$;
this homology group is called Schur multiplier of $G$.
\end{definition}

\begin{theorem}\label{schur}\cite{Ber}
If $(\Gamma,f,A)$ is a Schur cover of $G$, then every projective
representation $P$ of $G$ lifts to a linear representation of
$\Gamma$. Conversely, any linear representation of $\Gamma$ where A
acts by scalar matrices is a lift of a projective representation of
$G$.
\end{theorem}

\begin{remark}
Schur multiplier is an invariant of $G$ while the Schur cover is not
uniquely defined. But by last theorem, given a Schur cover of $G$,
all the projective representations of $G$ can be realized by linear
representations of $\Gamma$.
\end{remark}

\begin{definition}
Two projective representations of $G$ are called projective
equivalent if they are conjugated in $\mathbb{PGL}(n,\CC)$.
\end{definition}

Any projective representation of $G$ is given by a morphism of short
exact sequences:

\begin{equation}
\vcenter{\xymatrix{  0 \ar[r] & \CC^* \ar[r] & \mathbb{GL}(8,\CC)
\ar[r] &\mathbb{PGL}(8,\CC) \ar[r] & 1 \\
0\ar[r] & K \ar[r]\ar[u] & \Gamma \ar[r]\ar[u] & G \ar[r]\ar[u] & 1
}}
\end{equation}

where $\Gamma$ is a Schur cover of $G$. Usually the map $\tau$ is
not injective. Consider the short exact sequence:

\[
\begin{CD}
\ 0 @>>>K/Ker(\tau) @>>> \Gamma/Ker(\tau) @>>>G @>>>1
\end{CD}
\]
Here $K/Ker(\tau)$ is a cyclic group. By theorem \ref{schur},
projective representations of $G$ are in one to one correspondence
with linear representations of $\Gamma/Ker(\tau)$.

\begin{definition}
We say that a finite group $G\subset \mathbb{PGL}(8,\CC)$ has an
\emph{allowable action} if $G$ acts freely on some smooth complete
intersection $X$ of four quadrics in $\PP^7$. We will call the
correspondent $G$-action linear allowable action if $G$ can be
lifted to a subgroup of $\mathbb{GL}(8,\CC)$. Similarly, if the
variety $X$ is singular with ordinary double points, we say that $G$
has a \emph{semi-allowable action}.
\end{definition}

\begin{prop}\label{bk}
If $G$ has an allowable or semi-allowable action then $|G|$ divides
256.
\end{prop}

\begin{proof}
In [BKa], Browder and Katz proved a general theorem about free
action of finite groups on projective varieties:
\begin{theorem}\cite{BKa}
Let X be a projective variety in $\PP^n$ and $G$ is a finite
subgroup of $\mathbb{PGL}(n+1,\CC)$. If $G$ acts freely on X then,
$|G|$ divides the square of the degree of $X$.
\end{theorem}
We are considering complete intersections of four quadrics $X$ in
$\PP^7$, which have degree 16. By theorem of Browder and Katz, if
$G$ acts freely on $X$ then $|G|$ divides 256.
\end{proof}
\begin{remark}
Later we are going to argue the maximal order of $G$ is 64.
\end{remark}
If $G$ has an allowable action on $X$, then $H^0(X,\cO(1))$ becomes
a projective representation of $G$. We denote this vector space by
$V$. By theorem \ref{schur}, the group $\Gamma/Ker(\tau)$ acts
linearly on V with the cyclic subgroup $K/Ker(\tau)$ acting by
scalar matrices. By Holomorphic Lefschetz formula, those elements in
$\Gamma$ but not in $K/Ker(\tau)$ have trace zero. If we fix a
generator $\sigma$ of $K/Ker(\tau)$ of order $2^d$, then it should
act on $V$ as a scalar matrix $\xi I$ where $\xi$ is a primitive
$2^d$-th root of unity and $I$ stands for identity matrix. Let's
denote the trace vector of $\Gamma$ for a given representation $V$
by $t^{\Gamma}_{V}$. All the entries in $t^{\Gamma}_{V}$ are zero
except those corresponding to the conjugacy classes
$\{(\sigma^k),k=0,1,...,2^d-1\}$. These conjugacy classes have trace
$8\xi^{k}$. Similarly entries of $t^{\Gamma}_{H^0(X,\cO(2))}$ are
$32\xi^{2k}$ for conjugacy classes $\{(\sigma^k),k=0,1,...,2^d-1\}$
and zero otherwise. We can also compute the trace vector of the
induced representation $Sym^2(V)$ and denote it by
$t^{\Gamma}_{Sym^2(V)}$. The difference vector
$v=t^{\Gamma}_{Sym^2(V)}-t^{\Gamma}_{H^0(X,\cO(2))}$ is the trace
vector for the sub representation spanned by the four quadrics. The
assumption that $G$ acts freely on $X$ will force $t^{\Gamma}_{V}$
and $v$ to be group characters.
\begin{definition}
We say a central extension
\[
\begin{CD}
\ 0 @>>>K/Ker(\tau) @>>> \Gamma/Ker(\tau) @>>>G @>>>1
\end{CD}
\]
satisfies \emph{Lefschetz condition} if the trace vectors
$t^{\Gamma}_{V}$ and $v$ defined above are both group characters.
\end{definition}
\begin{prop}
If $G$ has a semi-allowable action then it satisfies Lefschetz
condition.
\end{prop}
\begin{proof}
Apply holomorphic Lefschetz formula to $\pi^*(\cO(1))$ and
$\pi^*(\cO(2))$ on the resolution $\pi: \widehat{X} \to X$.
\end{proof}
\begin{remark}
A priori, Lefschetz condition is only necessary but not sufficient
for $G$ to have allowable action. We still need to check the fixed
loci of $G$ in $\PP^7$ don't intersect with $X$ in order to verify
the freeness. However, in our cases it turns out that all the groups
satisfying Lefschetz condition are allowable when $|G|<64$.
When$|G|=64$ the necessity of Lefschetz condition follows from the
fact that ordinary double points are rational singularities. Details
are left to the readers.
\end{remark}

\section{Classification algorithm}\label{alg}
Our target is to classify the allowable and semi-allowable actions
on complete intersections of four quadrics in $\PP^7$ up to
projective equivalence. In this section we describe the scheme of
our algorithm. Some codes of the algorithms can be found in Appendix
of \cite{Hua}.

Recall that every projective representation gives a commutative diagram:\\

\begin{equation}
\vcenter{\xymatrix{  0 \ar[r] & \CC^* \ar[r] & \mathbb{GL}(8,\CC)
\ar[r] &\mathbb{PGL}(8,\CC) \ar[r] & 1 \\
0\ar[r] & K \ar[r]\ar[u]^\tau & \Gamma \ar[r]\ar[u] & G \ar[r]\ar[u]
& 1 }}
\end{equation}

where $\Gamma$ is a Schur cover of $G$ and $K$ is its Schur
multiplier. Generally $K$ is quite big but the exponent of $K$ is
controlled by order of $G$ by the following lemma.
\begin{lemma}
Let $G$ be a finite group and $K$ be its Schur multiplier. Denote
exponent of $K$ by $e$. Then $e^2$ divides $|G|$.
\end{lemma}
\begin{proof}
See \cite{Ber}.
\end{proof}

This lemma tells us the cyclic group $K/Ker(\tau)$ in the central
extension
\[
\begin{CD}
\ 0 @>>>K/Ker(\tau) @>>> \Gamma/Ker(\tau) @>>>G @>>>1
\end{CD}
\]
has order at most eight. By theorem \ref{schur}, given a group $G$
of order less or equal to 64, all projective representations of $G$
can be lift to a linear representation of $\Gamma/Ker(\tau)$.

Now we will describe the algorithm for $|G|=64$. Lower order groups
are handled similarly.
\begin{lemma}
If $|G|=64$ and $G$ acts freely on $X$ then $|K/Ker(\tau)|\geq 4$.
\end{lemma}
\begin{proof}
If $K/Ker(\tau)$ has order $2$ then the sheaf $\cO(2)$ must be $G$
linearizable, i.e. $dim(H^0(X,\cO(2)))$ must be divisible by 64. But
$H^0(X,\cO(2))$ has dimension $32$.
\end{proof}

Following this lemma, it suffices to consider projective
representations of a 64 group $G$ given by the following two types
of central extensions.
\begin{enumerate}
\item[1)] A group $H$ of order $256$ with a subgroup $\ZZ/4$ acting as diagonal matrix $\xi_8^2 I$;
\item[2)] A group $H$ of order $512$ with a subgroup $\ZZ/8$ acting as diagonal matrix $\xi_8 I$.
\end{enumerate}
Again $I$ represents the $8\times8$ identity matrix and $\xi_8$ is a
primitive $8$-th root of unity.

Now we can summarize our algorithm step by step.
\begin{enumerate}
\item[$\mathbf{Step\ I}$]: Check the Lefschetz condition for the
central extensions
\[
\begin{CD}
\ 0 @>>>\ZZ/4 @>>> H @>>>G @>>>1
\end{CD}
\]
and
\[
\begin{CD}
\ 0 @>>>\ZZ/8 @>>> H @>>>G @>>>1
\end{CD}
\]
Let $H$ go over all groups of order $256$ and $512$ and produce all
$G$ that satisfy Lefschetz condition. We use the GAP(\cite{GAP})
library of finite groups of small order. There are $56092$ different
groups of order $256$ and $10494213$ order $512$ groups.

\item[$\mathbf{Step\ II}$]: For each group $G$ that appears in Step I,
compute all possible extensions of $G$ of the form:
\[
\begin{CD}
\ 0 @>>>K/Ker(\tau) @>>> \Gamma/Ker(\tau) @>>>G @>>>1
\end{CD}
\] for a fixed Schur Cover $\Gamma$.
This can be done by computing kernels of all the group characters of
$K$. By theorem \ref{schur}, such extensions are in one to one
correspondence with nonequivalent projective representations.

\item[$\mathbf{Step\ III}$]: Check Lefschetz condition on extensions above and output those that satisfy it.

\item[$\mathbf{Step\ IV}$]: Check the fixed loci of the group actions obtained
above and show they don't intersect X.

\item[$\mathbf{Step\ V}$]: Check that the generic complete intersection
has at most ODP singularities in the semi-allowable case or is
smooth in the allowable case.
\end{enumerate}

The final output of the algorithm is a list of projective
representations of groups with allowable or semi-allowable actions.
The same group might appear on this list for several times with
different projective representations. Computer algebra system
involving in our algorithm are GAP (\cite{GAP}) and
MACAULAY(\cite{Macaulay}). The results of these calculations are
presented in the next section.
\section{Results}

In this section we present the results of the algorithm of the last
section.
\begin{remark}
Many group theoretic computation in this paper are done in GAP. It
has a small group library where all groups of given order less than
2000 are listed. For instance the quaternion group $H_8$ is
represented by $(8,4)$ in GAP library, where 8 for its order and 4
for its index in GAP library.
\end{remark}

There are 8 nontrivial groups of order less and equal to 8. We will
see all of them have allowable actions except the dihedral group
$D_8$. Further all the order 8 allowable action are linear.

There are 14(resp. 51) non-isomorphic 16-groups(resp. 32-groups). In
the following tables we list all the allowable groups by their
indices, together with the extension $\Gamma/Ker(\tau)$ representing
the correspondent projective representation. We also give number of
allowable actions up to projective equivalence.

When the order of the group is less than 64, the generic element of
the family with allowable action is a smooth complete intersection
of four quadrics in $\PP^7$. However this is no longer true for
64-groups.

There are 267 different groups of order 64. In these 267 groups
there are five groups that are semi-allowable.
\begin{table}[!h]\label{tab}
\tabcolsep 2mm \caption{(semi-)allowable action of order 2 to 64}
\begin{center}
\begin{tabular}{|r@{}lr@{}lr@{}l|}
\hline Groups&&Extension&&Schur Multiplier&
\\ \hline
$\ZZ/2$&   &$\ZZ/2$& &id group& \\
$\ZZ/4$&   &$\ZZ/4$& &id group& \\
$\ZZ/2\times \ZZ/2$&   &$\ZZ/2\times \ZZ/2$& &id group& \\
$\ZZ/8$&   &$\ZZ/8$& &id group& \\
$\ZZ/2\times\ZZ/4$&      &$\ZZ/2\times\ZZ/4$,(16,3)&   &$\ZZ/2$& \\
$(\ZZ/2)^3$&  &$(\ZZ/2)^3$&     &$(\ZZ/2)^3$& \\
$H_8$&   &$H_8$&   &id group& \\ \hline

(16,2)& &(64,18)& &$\ZZ/4$&  \\
(16,4)& &(32,14)& &$\ZZ/2$&  \\
(16,5)& &(32,5)& &$\ZZ/2$& \\
(16,10)& &(32,22)& &$(\ZZ/2)^3$& \\
(16,12)& &(32,29)& &$(\ZZ/2)^2$& \\ \hline

(32,2)& &(64,18),(64,23)& &$(\ZZ/2)^3$& \\
(32,3)& &(128,6)& &$\ZZ/4$& \\
(32,4)& &(64,28)& &$\ZZ/2$& \\
(32,5)&  &(64,4)& &$(\ZZ/2)^2$& \\
(32,13)& &(64,46)& &$\ZZ/2$& \\
(32,21)& &(128,462)& &$(\ZZ/2)^2\times \ZZ/4$& \\
(32,35)& &(64,182)& &$(\ZZ/2)^2$& \\
(32,47)& &(64,224)& &$(\ZZ/2)^5$& \\ \hline

(64,2)& &$(\ZZ/8)^2\rtimes \ZZ/8$& &$\ZZ/8$& \\
(64,3)& &(256,321)& &$\ZZ/4$& \\
(64,68)&  &(256,4235)&  &$\ZZ/2\times \ZZ/4$& \\
(64,72)&  &(256,4222),(256,4233)&  &$(\ZZ/2)^2\times\ZZ/4$& \\
(64,179)&  &(256,6447)&   &$\ZZ/4$& \\ \hline
\end{tabular}
\end{center}
\end{table}

\begin{remark}
We want to explore a little more about these five 64-groups because
it turns out the geometry of them are particularly interesting. The
group $(64,2)$ is the abelian group $\ZZ/8\times \ZZ/8$. Its Schur
cover is the Heisenberg group $(\ZZ/8)^2\ltimes \ZZ/8$. The group
$(64,3)$ is a semi-direct product of two copies of $\ZZ/8$ and
$(64,179)$ is a semi-direct product of quaternion group $H_8$ and
$\ZZ/8$. These first 3 groups all contain a maximal abelian subgroup
$\ZZ/4\times \ZZ/8$, which has GAP index $(32,3)$. It was observed
in \cite{BHua1} that these three 64-groups act on the same family.
This is a two dimensional subfamily of the three dimensional family
with $(32,3)$ action, which is invariant under certain involution.

The other two groups $(64,68)$ and $(64,72)$ don't have obvious
semi-direct product structures. Both of them contain a maximal
abelian subgroup $(\ZZ/4)^2\times \ZZ/2$, which has GAP index
$(32,21)$. These two groups act on a different two dimensional
family (See theorem \ref{p11}).
\end{remark}
\begin{remark}
All groups listed in Table \ref{tab} are subgroups of these five
64-groups with only two exceptions: $(32,4)$ and $(32,5)$. In
$(32,2)$ case, we are not sure whether both projective
representations are induced from representations of 64-groups. It
turns out all the actions for $|G|\leq 32$ in Table \ref{tab} are
allowable. When $|G|=64$, they are semi-allowable.
\end{remark}

\begin{remark}
The readers might observe that the 32-group $(32,2)$ and the
64-group $(64,72)$ have two different projective representations,
i.e. there are two non-conjugated embeddings of these finite groups
into $\mathbb{PGL}(8,\CC)$. Recall that projective representations
are one to one correspondent with central extensions. They are
quotient groups of some Schur Cover. It is a natural question to ask
that whether these two representation can be identified by some
outer automorphism of the group. It turns out that the two different
projective representations of $(64,72)$ are identified by some outer
automorphism of $(64,72)$. In other words, these are two different
ways of parameterizing the same subgroup of $\mathbb{PGL}(8,\CC)$
(see section 6).
\end{remark}
By proposition \ref{bk} the maximal order of allowable action we can
get is 256. Suppose there is an order 128 semi-allowable group. Then
all its 64 subgroups must be semi-allowable. By a GAP calculation we
check that there are no 128-groups, all of whose order 64 subgroups
are among $\{(64,2),(64,3),(64,68),(64,72),(64,179)\}$. Hence there
is no (semi-)allowable group of order bigger than 64.

\section{Complete intersection varieties}\label{ci}
In the last section we found five semi-allowable 64-groups. The
following two theorems show that three of them act freely on a two
dimensional family of complete intersections of four quadrics in
$\PP^7$, and the other two groups act freely on a different
dimension two family.

\begin{theorem}\label{p2}
Let $X$ be complete intersection of four quadrics:
$$
\begin{array}{l}
q_1= t_1(x_1^2+x_5^2) + t_2(x_2x_8+x_4x_6)+t_3 x_2x_7\\
q_2= t_1(x_2^2+x_6^2) + t_2(x_3x_1+x_5x_7)+t_3 x_3x_8\\
q_3= t_1(x_3^2+x_7^2) + t_2(x_4x_2+x_6x_8)+t_3 x_4x_1\\
q_4= t_1(x_4^2+x_8^2) + t_2(x_5x_3+x_7x_1)+t_3 x_5x_2.\\
\end{array}
$$
There are three groups $G_1,G_2,G_3$ contained in
$\mathbb{PGL}(8,\CC)$. Group $G_1$ is generated by $\tau$ and
$\sigma$ where $\sigma=(12345678)$ is permutation of the coordinates
$x_i$ and $\tau(x_i)=\xi^{i-1} x_i$ with $\xi$ a primitive $8$-th
root of unity. Group $G_2$ generated by $\tau$ and
$\sigma_1=(18325476)$. Then $G_2$ is a nonabelian group isomorphic
to a semidirect product of two copies of $\ZZ/8$. Group $G_3$ is
generated by $\tau$ and the permutations $\sigma_2= (1357)(2468)$
and $\sigma_3=(1256)(4387)$. It is a nonabelian group isomorphic to
a semidirect product of normal subgroup $\ZZ/8\ZZ$ generated by
$\tau$ and the quaternion group $H_8$ generated by $\sigma_2$ and
$\sigma_3$. As in remark 6.3, $G_1=(64,2)$, $G_2=(64,3)$ and
$G_3=(64,179)$. They act on $X$ without fixed points.
\end{theorem}
\begin{proof}
See \cite{BHua1}.
\end{proof}

Now we introduce the other two groups of order 64 having
semi-allowable actions. Define groups $G_4,G_5,{G_5}^{'}$ as
following subgroups of $\mathbb{GL}(8,\CC)$.
Group $G_4$ generated by coordinates transformations $\sigma_1,\sigma_2,\sigma_3$, where\\
$$\sigma_1: (x_1,\ldots,x_8)\mapsto(\xi x_7,\xi x_8,\xi^3 x_5,\xi^3 x_6,-\xi x_3,-\xi x_4,\xi^3 x_1,\xi^3
x_2)$$
$$\sigma_2: (x_1,\ldots,x_8)\mapsto(-x_2,i x_1,-x_4,-i x_3,-i x_6,x_5,i
x_8,x_7)$$
$$\sigma_3: (x_1,\ldots,x_8)\mapsto(\xi^3 x_5,-\xi^3 x_6,-\xi x_7,\xi x_8,\xi^3 x_1,-\xi^3 x_2,-\xi x_3,\xi
x_4)$$
Group $G_5$ is generated by $\sigma_3,\sigma_4,\sigma_5$ where\\
$$\sigma_4: (x_1,\ldots,x_8)\mapsto(\xi x_7,\xi x_8,-\xi^3 x_5,\xi^3 x_6,-\xi x_3,-\xi x_4,\xi^3 x_1,-\xi^3
x_2)$$
$$\sigma_5: (x_1,\ldots,x_8)\mapsto(\xi^3 x_6,\xi^3 x_5,-\xi x_8,\xi x_7,\xi^3 x_2,\xi^3 x_1,\xi x_4,-\xi
x_3)$$ Group ${G_5}^{'}$ is generated by
${\sigma_3},{\sigma_4},{\xi\sigma_5}$.

These three groups $G_4,G_5,{G_5}^{'}$ all have order 256. Their
indices in GAP are respectively (256,4235), (256,4222) and
(256,4233). The corresponding projective linear groups in
$\mathbb{PGL}(8,\CC)$ are $(64,68)$ and $(64,72)$. The last two
groups $G_5,{G_5}^{'}$ have the same projective group (64,72).
\begin{remark}
The two groups $G_5,{G_5}^{'}$ lead to two nonequivalent projective
representations of $(64,72)$(See \ref{tab}). However, the
corresponding projective subgroups of $\mathbb{PGL}(8,\CC)$ are the
same, i.e. these two representations only differ by an outer
automorphism of $(64,72)$.
\end{remark}

By abusing notations, from now on we denote the projectivizations by
$G_4$ and $G_5$.
\begin{theorem}\label{p11}
Let $X$ be a complete intersection of four quadrics in $\PP^7$ cut
out by:
$$
\begin{array}{l}\label{e11}
q_1= t_1(x_1^2+x_2^2)-t_2(x_3^2+x_4^2)+t_1(x_5^2+x_6^2)+t_2(x_7^2+x_8^2)\\
q_2= -t_2(x_1^2+x_2^2)+t_1(x_3^2+x_4^2)+t_2(x_5^2+x_6^2)+t_1(x_7^2+x_8^2)\\
q_3= s_1(x_1^2-x_2^2)-s_2(x_3^2-x_4^2)+s_1(x_5^2-x_6^2)+s_2(x_7^2-x_8^2)\\
q_4= -s_2(x_1^2-x_2^2)+s_1(x_3^2-x_4^2)+s_2(x_5^2-x_6^2)+s_1(x_7^2-x_8^2).\\
\end{array}
$$
The groups $G_4$ and $G_5$ introduced above act freely on $X$.
\end{theorem}

\begin{proof}
We only prove the theorem for $G_4$. The argument for $G_5$ is
completely analogous. Consider central extension
\[
\begin{CD}
\ 0 @>>>\ZZ/4 @>>>(256,4235)@>>> G_4 @>>>1\\
\end{CD}
\]

The group $(256,4235)$ has 46 irreducible representations, indexed
by $X_1...X_{46}$. In particular, $X_1,\ldots,X_{16}$ are one
dimensional irreducible representations, $X_{17},...,X_{44}$ are two
dimensional irreducible representations and $X_{45},X_{46}$ are
eight dimensional irreducible representations. We identify $V$ with
$H^0(X,\cO(1))$. Holomorphic Lefschetz formula force $V$ to be the
irreducible representation $X_{45}$. The second symmetric product of
$V$ has decomposition:
$$Sym^2(V)=(\oplus_{i\in I}
X_i)\oplus X_{35}^{\oplus 2}\oplus X_{36}^{\oplus 2}\ for\
I=\{19,20,21,22,25,26,27,28,33,34,41,42,43,44\}$$ The sub
representation spanned by the four quadrics has decomposition
$X_{35}\oplus X_{36}$, again follow from holomorphic Lefschetz
formula. Pick a basis $(x_1,\ldots,x_8)$ for $V=X_{45}$. We get an
induced basis for $Sym^2(V)$. They are homogenous quadratic
polynomials in $x_1,\ldots,x_8$. In particular,
$$X_{35}^{\oplus 2}=
Span\{x_1^2+x_2^2+x_5^2+x_6^2,x_3^2+x_4^2+x_7^2+x_8^2\} \oplus
Span\{x_7^2+x_8^2-x_3^2-x_4^2,x_5^2+x_6^2-x_1^2-x_2^2\}.$$
Respectively,
$$X_{36}^{\oplus 2}=
Span\{x_1^2-x_2^2+x_5^2-x_6^2,x_3^2-x_4^2+x_7^2-x_8^2\} \oplus
Span\{x_7^2-x_8^2-x_3^2+x_4^2,x_5^2-x_6^2-x_1^2+x_2^2\}.$$ These
polynomials give the cut out equations \ref{e11}. It is clear from
these equations that parameter space of this two dimensional family
is a subset of $\PP^1\times\PP^1$ where $(t_1:t_2)$ and $(s_1:s_2)$
are homogeneous coordinates of each $\PP^1$.

To show $G_4$ acts without fixed points, we need to check the
intersection of the fix loci of all conjugacy classes of $G_4$ with
$X$ are empty. It is easy to see this is the case for generic choice
of $t_1,t_2,s_1,s_2$.
\end{proof}

\begin{remark}
We have mentioned $(64,72)$ has two different projective
representations $(256,4222)$ and $(256,4233)$. A calculation shows
both of them act freely on this family.
\end{remark}

\section{Resolutions of Singularities}
We will investigate more about the geometry of these two families.
Let $X$ be a complete intersection of four quadrics cut out by
equations in theorem \ref{p2}. We have seen in last section three
64-groups $G_1$, $G_2$ and $G_3$ act freely on $X$. This family was
first discovered by Gross and Popescu. In \cite{GPo}, they studied
the birational geometry of $X$, including the resolution of
singularities. They have proved the following theorem in the case of
$G_1$.
\begin{theorem}
The singular Calabi-Yau 3-fold $X$ has an equivariant small
projective resolution $\widetilde{X}$, i.e. $\widetilde{X}$ is a
smooth projective Calabi-Yau 3-fold with free actions by $G_1$,
$G_2$and $G_3$. The resolution $\widetilde{X}$ has Hodge numbers
$h^{1,1}=2$ and $h^{1,2}=2$. Furthermore, $\widetilde{X}$ contains a
pencil of abelian surfaces with polarization $(1,8)$.
\end{theorem}

In this section we obtain a similar result for the family in theorem
\ref{p11}. We will prove the generic element $X$ in this family also
has an equivariant small projective
resolution. Recall $X$ is cut out by equations:\\
$$
\begin{array}{l}
q_1= t_1(x_1^2+x_2^2)-t_2(x_3^2+x_4^2)+t_1(x_5^2+x_6^2)+t_2(x_7^2+x_8^2)\\
q_2= -t_2(x_1^2+x_2^2)+t_1(x_3^2+x_4^2)+t_2(x_5^2+x_6^2)+t_1(x_7^2+x_8^2)\\
q_3= s_1(x_1^2-x_2^2)-s_2(x_3^2-x_4^2)+s_1(x_5^2-x_6^2)+s_2(x_7^2-x_8^2)\\
q_4= -s_2(x_1^2-x_2^2)+s_1(x_3^2-x_4^2)+s_2(x_5^2-x_6^2)+s_1(x_7^2-x_8^2).\\
\end{array}
$$
The jacobian matrix of it is
$$\left(
  \begin{array}{cccccccc}
   t_1x_1 & t_1x_2 & -t_2x_3 & -t_2x_4 & t_1x_5 & t_1x_6 & t_2x_7 & t_2x_8 \\
   -t_2x_1 & -t_2x_2 & t_1x_3 & t_1x_4 & t_2x_5 & t_2x_6 & t_1x_7 & t_1x_8 \\
   s_1x_1 & -s_1x_2 & -s_2x_3 & s_2x_4 & s_1x_5 & -s_1x_6 & s_2x_7 & -s_2x_8 \\
   -s_2x_1 & s_2x_2 & s_1x_3 & -s_1x_4 & s_2x_5 & -s_2x_6 & s_1x_7 & -s_1x_8 \\
  \end{array}
\right)$$

A point on $X$ is singular if and only if this matrix is
degenerated.

\begin{lemma}
A point $P\in X$ is singular if and only if exactly four coordinates
out of $(x_1,..,x_8)$ are zero.
\end{lemma}
\begin{proof}
Let $P=(x_1:\ldots:x_8)$ be a point on $X$. Observe that $P$ has at
most four zeros in coordinates because otherwise $P$ can't sit on
$X$ with generic choices of $t_i$ and $s_i$. We first prove if $P$
has 4 zeros then it must be a singular point. Let $\mu$ be a subset
of four distinct numbers in $\{1,\ldots,8\}$. Denote its complement
by $\bar{\mu}$. Let $P_{\mu}$ be a point with $\{x_i=0|i\in \mu\}$
and $J_{\bar{\mu}}$ be the four by four minor of the jacobian matrix
by picking the $\bar{\mu}$-th columns. Since all equations
$q_1,\ldots,q_4$ consist of square terms, the jacobian matrix $J$ is
equivalent with the coefficient matrix up to elementary
transformation. So $J_{\bar{\mu}}$ degenerates if and only if
$q_1,\ldots,q_4$ has nonzero solution of the form
$\{(x_1,\ldots,x_8)|x_i\neq 0\ for\ i\in \mu\ x_i=0\ for\ i\in
\bar{\mu}\}$. This proves the first direction.

If $P$ is a singular point, we pick a $\mu$ such that
$\{x_i\neq0|i\in \mu\}$. Since $P$ is singular the jacobian matrix J
evaluated at $P$ degenerates. In particular, $J_{\mu}$ degenerates.
If the other coordinates $\{x_i|i\in \bar{\mu}\}$ don't vanish
simultaneous, then we will get one parameter family of solutions
along which the jacobian matrix degenerates. However the singular
loci has dimension zero generally. So all the other coordinates must
be zero.
\end{proof}
Given $(1468)$, $(1367)$, $(1457)$, $(2467)$, $(2357)$, $(2458)$,
$(1358)$ and $(2368)$ as subsets of $\{1,\ldots,8\}$. The
corresponding points $P_{\mu}$ for $\mu$ equals any of these eight
sets are singular points of $X$. Since $X$ is cut out by degree two
equations, there are exactly 8 solutions for a given set $\mu$.
Hence each combinations give 8 singular points. These 64 points form
group orbits, for both $G_4$ and $G_5$. We will see later these 64
singularities are ordinary double points. Let's fix a set, say
$(1468)$. The corresponding singular points are
$(0:y_2:y_3:0:y_5:0:y_7:0)$. Plug $y_1=y_6=y_4=y_8=0$ into equations
\ref{e11}, we get:
$$
\begin{array}{l}
q_1= t_1y_2^2-t_2y_3^2+t_1y_5^2+t_2y_7^2=0\\
q_2= -t_2y_2^2+t_1y_3^2+t_2y_5^2+t_1y_7^2=0\\
q_3= -s_1y_2^2-s_2y_3^2+s_1y_5^2+s_2y_7^2=0\\
q_4= s_2y_2^2+s_1y_3^2+s_2y_5^2+s_1y_7^2=0.\\
\end{array}
$$

Solving $t_1,t_2,s_1,s_2$ by $y_i$, we rewrite the original
equations as:
$$
\begin{array}{lcr}
q_1&=& (y_3^2-y_7^2)(x_1^2+x_2^2)-(y_2^2+y_5^2)(x_3^2+x_4^2)\\&&+(y_3^2-y_7^2)(x_5^2+x_6^2)+(y_2^2+y_5^2)(x_7^2+x_8^2)\\
q_2&=& -(y_2^2+y_5^2)(x_1^2+x_2^2)+(y_3^2-y_7^2)(x_3^2+x_4^2)\\&&+(y_2^2+y_5^2)(x_5^2+x_6^2)+(y_3^2-y_7^2)(x_7^2+x_8^2)\\
q_3&=& (y_3^2-y_7^2)(x_1^2-x_2^2)-(y_5^2-y_2^2)(x_3^2-x_4^2)\\&&+(y_3^2-y_7^2)(x_5^2-x_6^2)+(y_5^2-y_2^2)(x_7^2-x_8^2)\\
q_4&=& -(y_5^2-y_2^2)(x_1^2-x_2^2)+(y_3^2-y_7^2)(x_3^2-x_4^2)\\&&+(y_5^2-y_2^2)(x_5^2-x_6^2)+(y_3^2-y_7^2)(x_7^2-x_8^2).\\
\end{array}
$$
Additionally $y_2,y_3,y_5,y_7$ satisfy a degree four relation
$y_3^4-y_7^4=y_2^4+y_5^4$.

These computations show that positions of the 64 singular points
uniquely determine the family of complete intersections.

No we will describe the explicit equivariant crepant resolution for
$X$ for $G_4$.
\begin{theorem}\label{res11}
There exist $G$-equivariant small resolutions
\xymatrix{\widetilde{X}\ar[r] & X} by blowing up a smooth
$G$-invariant abelian surface in $X$ for $G=G_4$ and $G=G_5$.
\end{theorem}
\begin{proof}
To construct such a small resolution, we need to find a Weil divisor
passing through the 64 ordinary double points, and invariant under
the action of $G_4$. Such a divisor is never Cartier since it is
locally cut out by more than one equation. By blowing up this
divisor we get the projective small resolution. Consider the
codimension one subscheme cut out by the following two equations.
$$
\begin{array}{l}
f_1=r_1x_1x_2-r_2x_3x_4+r_1x_5x_6+r_2x_7x_8\\
f_2=-r_2x_1x_2+r_1x_3x_4+r_2x_5x_6+r_1x_7x_8\\
\end{array}
$$
Notice equations $x_1x_2+x_5x_6$ and $x_3x_4+x_7x_8$ span the two
dimensional irreducible representation $X_{33}$ and $x_1x_2-x_5x_6$
and $x_3x_4-x_7x_8$ span the two dimensional irreducible
representation $X_{34}$. And $f_1, f_2$ are two generic elements in
$X_{33}\oplus X_{34}$. These two elements together with
$q_1,\ldots,q_4$ cut out an $G_4$ invariant surface in $X$. We
denote it by $S_{r_1,r_2}$. Under generic choices of coefficients
$r_1$ and $r_2$, this is a smooth abelian surface. If we pick any
one of $f_1$ and $f_2$ we will get unions of two abelian surfaces.
Hence $S_{r_1,r_2}$ is a Weil divisor but not Cartier. The abelian
surface $S_{r_1,r_2}$ has arithmetic genus $p_a$=-1, i.e. it is of
degree 16 in $\PP^7$. By varying $r_1$ and $r_2$ Any two such
surfaces intersect at the 64 singular points of $X$. It also follows
from the form of equations that they are ordinary double points. By
blowing up $S_{r_1,r_2}$, we get a smooth projective Calabi-Yau
threefold $\widetilde{X}$. Since $S_{r_1,r_2}$ is $G_4$-invariant
$\widetilde{X}$ also carries with a free $G_4$-action.
\end{proof}
\begin{remark}
In the case of $G_5$, we need to blow up a different Weil divisor
cut out by equations:
$$
\begin{array}{l}
f_1=r_1x_1x_5-r_2x_2x_6+r_1x_3x_7-r_2x_4x_8\\
f_2=-r_1x_1x_5+r_2x_2x_6+r_1x_3x_7-r_2x_4x_8\\
\end{array}
$$
Recall that there are two different allowable actions of $G_5$,
lifted to $G_5=(256,4222)$ and ${G_5}^{'}=(256,4233)$. Both of them
act on this surface, i.e. they have the same equivariant
resolutions.
\end{remark}

\begin{corollary}
The quotient variety $\widetilde{X}_{G_4}/G_4$(resp.
$\widetilde{X}_{G_5}/G_5$) is a smooth projective Calabi-Yau
threefold with fundamental group $G_4$(resp. $G_5$).
\end{corollary}

Similar to \cite{BHua1} and \cite{Gpo}, this family $X$ also carries
a fibration structure of abelian surfaces.
\begin{prop}
The equations $f_1,f_2$ form a sub linear system of dimension one of
$\cO(2)$ with 64 base points exactly at the 64 ordinary double
points.
\end{prop}
\begin{proof}
We need to show $\phi: x\mapsto (f_1(x) : f_2(x))$ is a rational map
defined outside the 64 ordinary double points. It is obvious $\phi$
is defined at $X\setminus S_{r_1,r_2}$. For any points on
$S_{r_1,r_2}$ that are not the 64 ordinary double points, $f_1$ and
$f_2$ have a common divisor, i.e. $S_{r_1,r_2}$ is cut out locally
just one equation. By dividing out the common divisor we extend
$\phi$ everywhere except the 64 ordinary double points.
\end{proof}
\begin{remark}
Consider the space of quadrics spanned by $q_1,\ldots,q_4$ together
with $f_1,f_2$. These equations cut out a $(2,4)$ polarized abelian
surfaces in $\PP^7$ (See \cite{Ba} for more about this abelian
surface). Any four linear independent equations of these six cut out
a Calabi-Yau complete intersection with 64 ordinary double points.
However only a two dimensional subfamily has free actions of $G_4$
and $G_5$.
\end{remark}

\begin{remark}\label{24}
Let $X$ be the Calabi-Yau threefold cut out by equations in theorem
\ref{p11}. It contains a pencil of $(2,4)$ polarized abelian
surfaces \cite{Ba}. Give the small resolution $\widetilde{X}$ in
theorem \ref{res11}. The Calabi-Yau threefold $\widetilde{X}$ has
Hodge number $h^{1,1}=10$ and $h^{1,2}=10$. As we stated in the last
remark, only a two dimensional subfamily in this ten dimensional
family has free actions of $G_4$ and $G_5$. A similar argument to
Remark $4.11$ in \cite{GPo} can be applied to compute the Hodge
number of the quotient variety $\widetilde{X}/G$. We expect the
quotient to have Hodge number $h^{1,1}=2,h^{1,2}=2$.
\end{remark}

\end{document}